\documentclass[12pt]{article}
\usepackage{amsfonts,amssymb,epsfig}
\usepackage{amsmath}

\date{}
\begin{document}
\newtheorem{df}{Definition}
\newtheorem{thm}{Theorem}
\newtheorem{lm}{Lemma}
\newtheorem{pr}{Proposition}
\newtheorem{co}{Corollary}
\newtheorem{re}{Remark}
\newtheorem{note}{Note}
\newtheorem{claim}{Claim}
\newtheorem{problem}{Problem}

\def\R{{\mathbb R}}

\def\E{\mathbb{E}}
\def\calF{{\cal F}}
\def\N{\mathbb{N}}
\def\calN{{\cal N}}
\def\calH{{\cal H}}
\def\n{\nu}
\def\a{\alpha}
\def\d{\delta}
\def\t{\theta}
\def\e{\varepsilon}
\def\t{\theta}
\def\pf{ \noindent {\bf Proof: \  }}
\def\trace{\rm trace}
\newcommand{\qed}{\hfill\vrule height6pt
width6pt depth0pt}
\def\endpf{\qed \medskip} \def\colon{{:}\;}
\setcounter{footnote}{0}

\def\Lip{{\rm Lip}}

\renewcommand{\qed}{\hfill\vrule height6pt  width6pt depth0pt}

\title{A quantitative version of the commutator theorem for zero trace matrices
}

\author{William B. Johnson\thanks{Supported in part by NSF DMS-1001321 and U.S.-Israel Binational Science Foundation
 }, Narutaka Ozawa\thanks{Supported in part by JSPS}, Gideon Schechtman\thanks{Supported in part by U.S.-Israel Binational Science Foundation. Participant NSF Workshop in Analysis and Probability, Texas A\&M University
 } } \maketitle

\begin{abstract}

Let $A$ be a $m\times m$ complex matrix with zero trace and let $\e>0$. Then there are $m\times m$ matrices $B$ and $C$ such that $A=[B,C]$ and $\|B\|\|C\|\le K_\e m^\e\|A\|$ where $K_\e$ depends only on $\e$. Moreover, the matrix $B$ can be taken to be normal.

\end{abstract}

\section{Introduction}
It is well known that a complex $m\times m$ matrix $A$ is a commutator (i.e., there are matrices $B$ and $C$ of the same dimensions as $A$ such that $A=[B,C]=BC-CB$) if and only if $A$ has zero trace. In such a situation clearly $\|A\|\le 2\|B\|\|C\|$ where $\|D\|$ denotes the norm of $D$ as an operator from $\ell_2^m$ to itself.

{\em Is it true that the converse holds?} That is, if $A$ has zero trace are there $m\times m$ matrices $B$ and $C$ such that $A=[B,C]$ and $\|B\|\|C\|\le K\|A\|$ for some absolute constant $K$?

Here we provide a weaker estimate: The above holds for $K=K_\e m^\e$ for every $\e>0$ where $K_\e$ depends only on $\e$. Moreover, the matrix $B$ can be taken to be normal.

The proof will be presented in the next section. It is self contained except for two facts. The first is a relatively easy result of Rosenblum \cite{ro} which gives a solution for $X$ of the matrix equation $A=SX-XT$ where all matrices are square and $S$ and $T$ have separated spectra in the sense that there is a domain $D$, whose boundary is a simple curve, which contains the spectrum of $S$ and is disjoint from the spectrum of $T$. The solution then is:
\[
X=\frac{1}{2\pi\imath}\int_{\partial D}(zI-S)^{-1}A(zI-T)^{-1}dz.
\]
The second fact is a heavy theorem of Bourgain and Tzafriri \cite{bt} related to restricted invertibility  of matrices and to the Kadison--Singer conjecture. It is stated as Theorem \ref{thm:bt} in the sequel.

The problem we discuss here was raised on MathOverFlow.net \cite{mo}. Although the MO discussion did not produce a solution to the problem, it did put the author in contact with one another and the discussion itself contains some useful tidbits.

\section{The main result}

Given $0<\e<1$, define a sequence of sets $\Lambda_n$ inductively: $\Lambda_1$ is the set of 4 points $\{\pm1\pm\imath1\}$ and
\[
\Lambda_n=\frac{1-\e}{2}\Lambda_{n-1}+\{\pm\frac{1+\e}{2}\pm\imath\frac{1+\e}{2}\}.
\]
Note that $\Lambda_n$ is a subset of the square $[-1,1]\times[-\imath,\imath]$ of cardinality $4^n$ and that it consists of a disjoint union of 4 sets each of which is a translate of $\frac{1-\e}{2}\Lambda_{n-1}$ and for each two of them their projection on either the real or imaginary axis is $2\e$ separated.

Given a $4^n\times4^n$ matrix $A$ with zero diagonal denote by $\mu(A)$ the smallest number $\mu$ such that there is a diagonal matrix $B$ with diagonal elements exactly the points of $\Lambda_n$ and a $4^n\times4^n$ matrix $C$ such that $A=[B,C]=BC-CB$ and $\|C\|\le\mu$. Note that since $A$ has zero diagonal, for each diagonal matrix $B$ with distinct diagonal entries $\{b_i\}$ such a matrix $C$ exist and its non diagonal entries are uniquely defined by $c_{ij}=a_{ij}/(b_i-b_j)$. Put also $\mu(4^n)=\max\mu(A)$ where the $\max$ ranges over all zero diagonal $4^n\times4^n$ matrices of norm one.

Similarly, for $m$ not necessarily of the form $4^n$, we denote by $\lambda(A)$ the smallest number $\lambda$ such that there is a diagonal matrix $B$ with diagonal elements in $[-1,1]\times[-\imath,\imath]$ and a $m\times m$ matrix $C$ such that $A=[B,C]=BC-CB$ and $\|C\|\le\lambda$. Put $\lambda(m)=\max\lambda(A)$ where the $\max$ ranges over all $m\times m$ matrices of zero diagonal and norm one.

Given a $m\times m$, $m=4^n$, matrix $A$ write it as a $4\times 4$ block matrix with blocks of size $4^{n-1}\times4^{n-1}$
\[
\left(
  \begin{array}{cccc}
    A_{11} & A_{12} & A_{13} & A_{14} \\
    A_{21} & A_{22} & A_{23} & A_{24} \\
    A_{31} & A_{32} & A_{33} & A_{34} \\
    A_{41} & A_{42} & A_{43} & A_{44} \\
  \end{array}
\right)
\]
\begin{claim}\label{claim:mu}
\[\mu(A)\le \frac{2}{1-\e}\max_{1\le i\le 4}\mu(A_{ii})+\frac{6\|A\|}{\e^2}.
\]
In particular
\[
\mu(4^n)\le \frac{2}{1-\e}\mu(4^{n-1})+\frac{6}{\e^2}.
\]
Also,
\begin{equation}\label{eq:1}
\lambda(A)\le \frac{2}{1-\e}\max_{1\le i\le 4}\lambda(A_{ii})+\frac{6\|A\|}{\e^2}\ \ \mbox{and}\ \ \lambda(4^n)\le \frac{2}{1-\e}\lambda(4^{n-1})+\frac{6}{\e^2}.
\end{equation}
\end{claim}
\pf
Let $B_{ii}$ be diagonal matrices with diagonal entries in $\Lambda_{n-1}$ and $C_{ii}$ $4^{n-1}\times4^{n-1}$ matrices with $A_{ii}=[B_{ii},C_{ii}]$ and $\|C_{ii}\|=\mu(A_{ii})$.
Let
\[\{B_{ii}^\prime\}_{i=1}^4=\left\{\frac{1-\e}{2}B_{\frac{a+1}2+b+2,
\frac{a+1}2+b+2}+\left(a\frac{1+\e}{2}+\imath b\frac{1+\e}{2}\right)I_{4^{n-1}}\right\}_{a,b=\pm1}
\]
(the order doesn't matter), and, for $i\not=j$, let $C^\prime_{ij}$ be defined (uniquely) by
\[
A_{ij}=B_{ii}^\prime C_{ij}^\prime-C_{ij}^\prime B_{jj}^\prime.
\]
Then by the result mentioned in the Introduction (see\cite{ro} or \cite{ma}),
\[
C_{ij}^\prime=\frac{1}{2\pi\imath}\int_{\partial D_{ij}}(zI-B_{ii}^\prime)^{-1}A_{i,j}(zI-B_{jj}^\prime)^{-1}dz
\]
where $D_{ij}$ is the boundary curve of any domain containing the spectrum of $B_{ii}^\prime$ and disjoint from the spectrum of $B_{jj}^\prime$. Since we can easily find such a curve of distance at least $\e$ from the spectra of $B_{ii}^\prime$ and $B_{jj}^\prime$ and of length $4+4\e<8$ we get that
$\|C_{ij}^\prime\|<\frac{2}{\e^2}\|A_{ij}\|$.

Let
$C_{ii}^\prime=\frac{2}{1-\e}C_{ii}$ and set
\[
B=\left(
    \begin{array}{cccc}
      B_{11}^\prime & 0 & 0 & 0 \\
      0 & B_{22}^\prime & 0 & 0 \\
      0 & 0 & B_{33}^\prime & 0 \\
      0 & 0 & 0 & B_{44}^\prime \\
    \end{array}
  \right)
  \]
  and
  \[C=(C_{ij}^\prime)_{i,j=1,2,3,4}.
  \]
  Then
  \[\|C\|\le\frac{2}{1-\e}\max_{i,i}\mu(A_{ii})+\frac{6}{\e^2}\|A\|.
  \]
This gives the claim for $\mu$  and the proof for $\lambda$ is almost identical.
\endpf

  In the proof of the main theorem we shall use the parameter $\lambda$. The reason we also included $\mu$ here is that the matrices $B$ in the proof for the property of $\mu$ depend only on $\e$ and not on the matrices $A$. Optimizing over $\e$ we get
  \begin{co}\label{co:1}
  (i) For each $m$ there is a $m\times m$ diagonal matrix $B$ with spectrum in the square $[-1,1]\times[-\imath,\imath]$ such that for each $m\times m$ matrix $A$ with diagonal zero there is a $m\times m$ matrix $C$ with norm at most $O((\log m)^3\sqrt m)\|A\|$ such that $A=[B,C]$.

  (ii) For each $m=4^n$ there is a subset $\Lambda_m$ of $[-1,1]\times[-\imath,\imath]$ such that any trace zero $m\times m$ matrix $A$ there is a normal matrix $B$ with spectrum $\Lambda_m$ and a matrix $C$ with norm at most $O((\log m)^3\sqrt m)\|A\|$ such that $A=[B,C]$.
\end{co}

\pf For each $0<\e<1$, $m$ of the form $4^n$, and an $m\times m$ matrix $A$ with norm 1 and zero diagonal, Claim \ref{claim:mu} gives, as long as $\frac{6}{\e^2}\le \frac{2\e}{1-\e}\mu(m/4)$, that
\[
\mu(m)\le 2\frac{1+\e}{1-\e}\mu(m/4).
\]
Let $k$ be the largest natural number smaller than $\log_4 m$ such that $\frac{6}{\e^2}\le \frac{2\e}{1-\e}\mu(m/4^k)$. (If no such $k$ exists take $k=\log_4 m$ and change the argument below a bit, getting a better estimate.) Then
\begin{eqnarray*}
\mu(m)&\le&(2\frac{1+\e}{1-\e})^k\mu(m4^{-k})\le(2\frac{1+\e}{1-\e})^k(\frac{2}{1-\e}\mu(m4^{-(k+1)})+\frac{6}{\e^2})\cr
&\le&(2\frac{1+\e}{1-\e})^k(\frac{6}{\e^3}+\frac{6}{\e^2})\le \frac{12}{\e^3}(2\frac{1+\e}{1-\e})^k.
\end{eqnarray*}
For $\e=\frac{1}{k}$ we get
\[
\mu(m)\le 12k^3 2^k(1+\frac{3}{k})^k.
\]
Since $k$ is at most $\log_4 m$ we get (i) to get (ii) use the fact (see e.g. \cite{fi} or \cite{ha}) that any trace zero matrix is unitarily equivalent to a matrix with zero diagonal.
\endpf

\begin{re}{\rm
The power ${1/2}$ of $m$ in the first part of Corollary \ref{co:1} can't be lowered. Indeed, if $B$ is any $m\times m$ diagonal matrix with spectrum in $[-1,1]\times [-\imath,\imath]$ then there are $i\not= j$ in $\{1,2,\cdots,m\}$ with $|i-j|\le \sqrt{8/m}$. If $A$ is the  $m\times m$ matrix with 1 in the $i,j$ place and zero elsewhere and $A=[B,C]$, then it is easy to see that the absolute value of the $i,j$ entry of $C$ is at least $\sqrt{m/8}$.}
\end{re}

Note that the constant $\frac{2}{1-\e}$ in (\ref{eq:1}) is what leads to the power $1/2$ of $m$ in the Corollary above. If we could replace it with $\frac{1}{1-\e}$ we could eliminate  the power of $m$ altogether and be left with only a $\log$ factor. The next Claim is a step in this direction. The Claim, which has a proof similar to the previous one, shows that if a zero diagonal $2m\times 2m$ matrix $A$ has its two $m\times m$ central submatrices having substantially different $\lambda$ values and the smaller one is substantially larger than the norm of the matrix, then $\lambda(A)$ is, up to a multiplicative constant close to $1$, basically the same as the larger of these two values. This will be used in the proof of the main theorem.

\begin{claim}\label{claim:2} Let
\[
A=\left(
    \begin{array}{cc}
      A_{11}& A_{12}\\
      A_{21}& A_{22}\\
    \end{array}
  \right)
  \]
be a $2m\times 2m$ matrix with zero diagonal where the $A_{ij}$ are all $m\times m$ matrices. Assume also that $\lambda(A_{ii})\le c_i$ where $c_1/c_2<1/4$. Then
\[
\lambda(A)\le (1+K((c_1/c_2)^{1/2}+\|A\|/c_1))c_2
\]
For some absolute constant $K>0$.
\end{claim}
\pf
Write $A_{ii}=B_{ii}C_{ii}-C_{ii}B_{ii}$, $i=1,2$ where the $B_{ii}$ are diagonal matrices with spectrum in $[-1,1]\times[-\imath,\imath]$ and $\|C_{ii}\|=\lambda(A_{ii})\le c_i$. Assume also that $c_1<c_2$. For any $1/2>\delta\ge c_1/c_2$ put
\[
B_{11}^\prime=(-1+\delta)I+\delta B_{11}, \ \ \ \ B_{22}^\prime=2\delta I+(1-2\delta)B_{22}
\]
and
\[
C_{11}^\prime=\delta^{-1}C_{11}, \ \ \ \ C_{22}^\prime=(1-2\delta)^{-1}C_{22}.
\]
Then $A_{ii}=B_{ii}^\prime C_{ii}^\prime-C_{ii}^\prime B_{ii}^\prime$ and the $B_{ii}^\prime$-s are diagonal matrices with spectrum in $[-1,1]\times[-\imath,\imath]$. Moreover, the spectrum of $B_{11}^\prime$ lies to the left of the vertical line $\Re z=-1 +2\delta$ and that of $B_{22}^\prime$ to the right of the vertical line $\Re z=-1+4\delta$. Also
\[
\max_{i=1,2}\|C_{ii}^\prime\|\le\max\{\delta^{-1}c_1,(1-2\delta)^{-1}c_2\}=\frac{c_2}{1-2\delta}.
\]
Define $C_{ij}^\prime$, $i\not= j\in {1,2}$, by
\[
A_{ij}=B_{ii}^\prime C_{ij}^\prime-C_{ij}^\prime B_{jj}^\prime
\]
then, by the same argument as in the proof of Claim \ref{claim:mu}, using Rosenblum's result,
$\|C_{ij}\|\le K\|A\|/\delta^2$ for some universal $K$.
Define
\[
B^\prime=\left(
    \begin{array}{cc}
      B^\prime_{11}& 0\\
    0& B^\prime_{22}\\
    \end{array}
  \right)
  \ \ \ and \ \ \
  C^\prime=\left(
    \begin{array}{cc}
      C^\prime_{11}& C^\prime_{12}\\
      C^\prime_{21}& C^\prime_{22}\\
    \end{array}
  \right)
\]
then $A=B^\prime C^\prime-C^\prime B^\prime$, $B$ is a diagonal matrix with spectrum in
$[-1,1]\times[-\imath,\imath]$ and
\[
\|C^\prime\|\le \frac{c_2}{1-2\delta}+\frac{K\|A\|}{\delta^2}.
\]
Taking $\delta= (c_1/c_2)^{1/2}$ we get that
\[
\lambda(A)\le (1+K((c_1/c_2)^{1/2}+\|A\|/c_1))c_2
\]
for some absolute constant $K$ (which, a careful examination of the proof shows, can be taken to be $4/\pi$).
\endpf

We next recall a theorem of Bourgain and Tzafriri \cite{bt}.
\begin{thm}\label{thm:bt} \cite{bt}. For some absolute constant $K>0$, if $A$ is a $m\times m$ matrix with zero diagonal then for all $\e>0$ there is a central (i.e., whose diagonal is a subset of the diagonal of $A$) submatrix $A^\prime$ of dimension $\lfloor \e^2m\times \e^2m\rfloor $ whose norm is at most $K\e \|A\|$. \hfil\break
Consequently, If $A$ is a norm one $2\cdot4^n\times 2\cdot4^n$ matrix with zero diagonal then for all $l\le n$ there are $4^l$ disjoint subsets $\sigma_i$ of $1,2,\dots,2\cdot4^n$ each of size $4^{n-l}$ such that all the submatrices corresponding to the entries in $\sigma_i\times \sigma_i$ have norm at most $K2^{-l}$.
\end{thm}

\begin{thm}\label{thm:main}
(i) For each $\e>0$ there is a constant $K_\e$ such that for all $m$
\[
\lambda(m)\le K_\e m^\e.
\]
(ii) For each $\e>0$ there is a constant $K_\e$ such that for all $m$ and every $m\times m$ zero trace matrix $A$ there is a normal matrix $B$ with spectrum in $[-1,1]\times [-\imath,\imath]$ and a matrix $C$ with norm at most $K_\e m^\e\|A\|$ such that $A=[B,C]$.
\end{thm}

\pf Let $A$ be a $2\cdot4^n\times 2\cdot4^n$ matrix with zero diagonal and norm one. Let $1\le l\le n$ and let $A^\prime$ be the $4^n\times 4^n$ submatrix corresponding to the entries in $\cup_{i=1}^{4^l}\sigma_i\times \cup_{i=1}^{4^l}\sigma_i$ where $\sigma_i$ are given by Theorem \ref{thm:bt}.
Let $A_{ii}^l$ denote the submatrix corresponding to the entries in $\sigma_i\times \sigma_i$, $i=1,2,\dots,4^l$. Divide $1,2,\dots,4^l$ into $4^{l-1}$ disjoint sets each a union of $4$ $\sigma_i$-s and let $A_{ii}^{l-1}$, $i=1,2,\dots,4^{l-1}$,  denote the $4^{n-l+1}\times 4^{n-l+1}$ submatrices corresponding to the entries corresponding to these sets. Continue in this manner to define $A_{ii}^s$, $i=1,2,\dots,4^s$ for each $s=0,1,2,\dots,l$ where for $s\ge 1$ $A_{ii}^s$ is a $4^{n-s}\times 4^{n-s}$ submatrix of one of the $A_{jj}^{s-1}$. Note that $A^\prime=A_{11}^0$.

Now, By Claim \ref{claim:mu} for each $\e>0$,
\begin{eqnarray*}
\lambda(A^\prime)&\le& \frac{2}{1-\e}\max_{1\le i\le 4}\lambda(A_{ii}^1)+\frac{6}{\e^2}\cr
&\le& \left(\frac{2}{1-\e}\right)^2\max_{1\le i\le 16}\lambda(A_{ii}^2)+\left(\frac{2}{1-\e}+1\right)\frac{6}{\e^2}\cr
&\le& \ \ \ \dots\dots\cr
&\le& \left(\frac{2}{1-\e}\right)^{l-1}\max_{1\le i\le 4^{l-1}}\lambda(A_{ii}^{l-1})+\left(\left(\frac{2}{1-\e}\right)^{l-2}+\cdots+ \frac{2}{1-\e}+1\right)\frac{6}{\e^2}\cr
&\le& \left(\frac{2}{1-\e}\right)^{l}\lambda(4^{n-l})K2^{-l}+\left(\left(\frac{2}{1-\e}\right)^{l-1}+\cdots+ \frac{2}{1-\e}+1\right)\frac{6}{\e^2},
\end{eqnarray*}
where the last step is the place we use Theorem \ref{thm:bt}. Now use Corollary \ref{co:1} to get that for some absolute constants $K$ (not necessarily the same in each row)
\begin{align}\label{eq:lambda}
\lambda(A^\prime)\le & K\left(\frac{1}{1-\e}\right)^{l}\lambda(4^{n-l})+l\left(\frac{2}{1-\e}\right)^{l-1}\frac{6}{\e^2}\\
\le & K\left(\frac{1}{1-\e}\right)^{l}(n-l)^32^{n-l}+l\left(\frac{2}{1-\e}\right)^{l-1}\frac{6}{\e^2}\nonumber.
\end{align}
For $\e=1/l$ we get
\[
\lambda(A^\prime)\le K((n-l)^32^{n-l}+l^32^l)
\]
and taking $l=n/2$ gives
\begin{equation}\label{eq:aprime}
\lambda(A^\prime)\le Kn^32^{n/2}=K(\log m)^3m^{1/4}.
\end{equation}
We managed to reduce the power of $m$ in the bound on $\lambda(A)$ from $m^{1/2}$ to $m^{1/4}$ but only for a large submatrix. Next we are going to utilize Claim \ref{claim:2} to get a similar bound for the whole matrix. Let $\sigma^c=\{1,2,\cdots,2\cdot4^n\}\setminus \cup_{i=1}^{4^l}\sigma_i$ and let $A^{\prime\prime}$ be the submatrix of $A$ with entries in
$\sigma^c\times \sigma^c$. Put $c_1=K(\log m)^3m^{1/4}$ and $c_2=\max\{K(\log m)^7m^{1/4},\lambda(A^{\prime\prime})\}$. Then $A$, $A_{11}=A^{\prime}$ and $A_{22}=A^{\prime\prime}$ satisfy the assumptions of Claim \ref{claim:2} with $c_1,c_2$. Consequently,
\[
\lambda(A)\le (1+K(\log m)^{-2})\max\{K(\log m)^7m^{1/4},\lambda(A^{\prime\prime})\}
\]
where we continue to use $K$ to denote a universal constant, possibly different in different occurrences, and for $m=4^n$, $n\ge 1$,
\[
\lambda(2m)\le (1+K(\log m)^{-2})\max\{K(\log m)^7m^{1/4},\lambda(m)\}.
\]
Repeating the argument again reducing from matrices of size $4^{n+1}\times 4^{n+1}$ to ones of size $2\cdot4^{n}\times 2\cdot4^{n}$ and combining with the above we get, for $m=4^n$,
\[
\lambda(4m)\le (1+K(\log m)^{-2})\max\{K(\log m)^7m^{1/4},\lambda(m)\}.
\]
Let $k\le m $ be the largest power of $4$ such that $\lambda(k)\le K(\log_4k)^7k^{1/4}$. Then
\[
\lambda(4m)\le\left(\prod_{s=\log_4k+1}^{\log_4m}(1+Ks^{-2})\right)K(\log k)^7k^{1/4}.
\]
For some other absolute constant $K$ this last quantity is at most $K(\log m)^7m^{1/4}$.
We thus improved the previous bound on $\lambda(m)$ (for $m=4^n$) to
\[
\lambda(m)\le K(\log m)^7m^{1/4}
\]
for some absolute $K$.

Repeating the argument one can improve the bound further: Go back to (\ref{eq:lambda}) and plug this new bound to get
\[
\lambda(A^\prime)\le
K\left(\frac{1}{1-\e}\right)^{l}(n-l)^72^{(n-l)/2}+l\left(\frac{2}{1-\e}\right)^{l-1}\frac{6}{\e^2}.
\]
For $\e=1/l$ we get
\[
\lambda(A^\prime)\le K((n-l)^72^{(n-l)/2}+l^32^l)
\]
and taking $l=n/3$ gives
\[
\lambda(A^\prime)\le Kn^72^{n/3}=K(\log m)^7m^{1/6}.
\]
replacing (\ref{eq:aprime}) with this new estimate and following the rest of the argument above leads to
\[
\lambda(m)\le K(\log m)^{11}m^{1/6}.
\]
Iterating, this leads to a bounds of the form:
\begin{equation}\label{eq:Kk}
\lambda(m)\le K_k(\log m)^{4k-1}m^{1/2k}
\end{equation}
for every $m=4^n$ and every positive integer $k$, where $K_k$ depends only on $k$.
This gives the statement of the theorem for $m$ being a power of $4$. For a general $m\times m$ zero diagonal matrix $A$, complete it to a $4^n\times 4^n$ matrix $A^\prime$ where $4^{n-1}<m\le 4^n$ by adding zero entries and keeping $A$ supported on $\{1,2,\cdots,m\}\times \{1,2,\cdots,m\}$. Apply the theorem to $A^\prime$ and note that the fact that $B$ is diagonal implies that we can assume that $C$ has non zero entries only in $\{1,2,\cdots,m\}\times \{1,2,\cdots,m\}$. This proves the first part of the theorem.
The second follows from the fact that any trace zero matrix is unitarily equivalent to a zero diagonal matrix.
\endpf

\section{Concluding remarks}

1. Recall that the paving conjecture states that for every $\e>0$ there is a positive integer $n(\e)$ such that any norm one zero diagonal matrix has a paving of length at most $n(\e)$ and norm at most $\e$. By a paving of $A$ we mean a block diagonal submatrix of $A$ whose diagonal is the same as that of $A$. The length of a paving is the number of blocks. Anderson \cite{an} showed that this conjecture is equivalent to the Kadison--Singer conjecture \cite{ks} on the extension of pure states. For s recent expository paper on these conjectures see \cite{ce}. 

It is clear from the proof above that if the paving conjecture holds with the right parameters than the proof can be simplified and the main result strengthened to get a polylog estimate on $\lambda(m)$. We next show that the reverse holds in a very strong sense. In particular if $\lambda(m)$ is bounded independently of $m$ then the paving conjecture holds.
\begin{claim}
Assume $A=[B,C]$ with $B$ a $m\times m$ diagonal matrix with spectrum in $[-1,1]\times [-\imath,\imath]$ and $C$ an $m\times m$ matrix. Then for every $0<\e<1$ $A$ has a
paving of length $\lfloor \frac{2}{\e}\rfloor^2$ and norm $\sqrt2\e\|C\|$.
\end{claim}
\pf
Partition $[-1,1]$ into $\lfloor \frac{2}{\e}\rfloor$ disjoint intervals $I_i$ of length at most $\e$ each. Let $B(i,j)$ be the central (diagonal) submatrix of $B$ whose diagonal entries are in
$I_i\times\imath I_j$, let $A(i,j)$ and $C(i,j)$ be the central submatrices of $A$ and $C$ respectively with the same support as $B(i,j)$. $A(i,j)$, \ $i,j=1,2,\cdots$, $\lfloor \frac{2}{\e}\rfloor$. $A(i,j)$, $i,j=1,2,\cdots, \lfloor \frac{2}{\e}\rfloor$, is a paving of $A$ and it is enough to prove that $\|A(i,j)\|\le \sqrt 2 \e\|C\|$.

Clearly $A(i,j)=[B(i,j),C(i,j)]$. Pick $i,j$, let $b$ be the center of the square $I_i\times\imath I_j$ and note that $bI-B(i,j)$ (with I the identity matrix of the same dimensions as $B(i,j)$) is a diagonal matrix with entries of absolute value at most $\e/\sqrt 2$. Therefore
\[
\|A(i,j)\|=\|(B(i,j)-bI)C(i,j)-C(i,j)(B(i,j)-bI)\|\le \sqrt 2 \e\|C\|.
\]
\endpf

\medskip

\noindent
2. A more careful examination of the proof of Theorem \ref{thm:main} shows that the constant we get in (\ref{eq:Kk}) is 
\[
\lambda(m)\le K^k(\log m)^{4k-1}m^{1/2k}
\]
for some absolute constant $K$.
Optimizing over $k$ gives
\[
\lambda(m)\le m^{K(\log\log m/\log m)^{1/2}}
\]
for some absolute $K$.

\medskip

\noindent
3. Although the problem we discuss seems basic enough not to need further motivation, we would like to indicate one. If any trace zero matrix $A$ could be written as $A=[B,C]$ with $\|B\|\|C\|\le K\|A\|$ for a universal $K$, then we would get a simple characterization of the commutators in an important class of $II_1$ factors, the Wright factors; an element there would be a commutator if and only if it has zero trace. See \cite{ds} for this and related matters.

%
%

\begin{tabular}{ll}
W.B. Johnson&N. Ozawa\\
Department of Mathematics&Research Institute for Mathematical Sciences\\
Texas A\&M University&Kyoto University\\
College Station, TX  77843 U.S.A.&Kyoto 606-8502, Japan\\
{\tt johnson@math.tamu.edu}&{\tt narutaka@kurims.kyoto-u.ac.jp}
\\
\end{tabular}

\bigskip

\begin{tabular}{l}
G. Schechtman\\
Department of Mathematics\\
Weizmann Institute of Science\\
Rehovot, Israel\\
{\tt gideon@weizmann.ac.il}\\
\end{tabular}

\end{document}